\documentclass[12pt]{article}
\usepackage{graphics}
\usepackage{amsmath}
\usepackage[latin1]{inputenc}
\usepackage{amssymb}
\usepackage{graphicx}
\usepackage{makeidx}
\usepackage{mathrsfs}
\usepackage{amsfonts}
\usepackage[mathscr]{eucal}
\usepackage{amsmath}
\def\pni{\par \noindent}
\def\vsh{\vskip 0.25truecm\noindent}
\def\vshm{\vskip 0.15truecm}
\def\vs{\vskip 0.5truecm}

\def\vsp{\vsh\pni}%%%%%%%%%%%%%%%

\def\vsn{\vshm\pni}
\def\cen{\centerline}

\def\sgn{\hbox{sign}\,}

\def\e{\hbox{e}}
\def\exp{\hbox{exp}}
\def\ds{\displaystyle}

\def\q{\quad}

\def\lt{\left} \def\rt{\right}

\def\d{\partial}
   \def\dt{\partial t}

\def\NN{{\rm I\hskip-2pt N}}

\def\RR{\vbox {\hbox to 8.9pt {I\hskip-2.1pt R\hfil}}\;}

\def\exp{{\rm exp}\,} \def\e{{\rm e}}
\def\L{{\cal L}} %%% Laplace Transform !!!!
\def\F{{\cal F}} %%% Fourier Transform !!!!
\begin{document}

\cen{{\bf FRACALMO PRE-PRINT: \ http://www.fracalmo.org}}
% \vsh
\cen{{\bf The European Physical Journal, Special Topics, Vol. 193  (2011) 119--132}}
% \vsh
\cen{{\bf Special issue:   Perspectives on Fractional Dynamics and Control}}
%% \vsh
\cen{{\bf Guest Editors: Changpin LI  and Francesco MAINARDI }}
\vsh
\hrule
% \end{center}
%%%%%%%%%%%%%%%%%%%%%%%%%%%%%%%%%%%%%%%%%%%%%%%%%%%%%%%%%%%%%%%%%%%%%%%%%
%%BEGINNING OF TEXT

%%%%%%%%%%%%%%%%%%%%%%%%%%%%%%%%%%%%%%%%%%%%%%%%%%%%%%%%%%%%%%%%%%%%%%%%%
\vskip 0.50truecm
\font\title=cmbx12 scaled\magstep2
\font\bfs=cmbx12 scaled\magstep1
\font\little=cmr10
\begin{center}
%%{\bfs Time-fractional derivatives in relaxation processes:}
%%\vs{\bfs a tutorial survey}\vvs
%%%%%%%%%%%%%%%%%%%%%%%%%%%%%%%%%%%%%%%%%%%%%%%%%%%%%%%%%%%%%%%%%%%%%%%%%
\cen{{\title Subordination Pathways to Fractional Diffusion}\footnote{This paper would be referred to as
{\it Eur. Phys. J. Special Topics} {\bf 193}, 119--132 (2011).}}
\vs
{Rudolf GORENFLO} $^{(1)}$ and
{Francesco MAINARDI}$^{(2)}$
%%%%%%%%%%%%%%%%
\vsh
$\null^{(1)}$
 {\little Department of Mathematics and Informatics, Free University  Berlin,} \\
{\little Arnimallee 3, D-14195 Berlin, Germany} \\
%%%%%%%%%%%%%%%% Phone: +49 30 4504 2247 \\
{\little E-mail: {\tt gorenflo@mi.fu-berlin.de}}
\\ [0.25 truecm]
$\null^{(2)}$
 {\little Department of Physics, University of Bologna, and INFN,} \\
              %%%  Sezione di Bologna, \\
{\little Via Irnerio 46, I-40126 Bologna, Italy} \\
%% Tel: +39-051.2091098 $\;$ Fax: +39-051.247244 $\;$\\
{\little Corresponding Author. E-mail: {\tt francesco.mainardi@unibo.it}} 
\\ [0.25 truecm]

\end{center}

\begin{abstract}
\noindent
The uncoupled  Continuous Time Random Walk ($CTRW$) in one space-dimension and under power law regime is splitted 
  into three distinct random walks:
  ($rw_1$), a random walk along the line of natural time, happening in operational time; 
  ($rw_2$), a random walk along the line of space, happening in operational time;
  ($rw_3$), the inversion of ($rw_1$), namely a random walk along the line of operational 
  time, happening in natural time.
  Via the general integral equation of $CTRW$ and appropriate rescaling, the transition to the diffusion limit
  is carried out for each of these three random walks. Combining the limits of ($rw_1$) and ($rw_2$) 
  we get the method of parametric
   subordination for generating particle paths, 
   whereas combination of ($rw_2$) and ($rw_3$) yields the subordination integral for 
   the sojourn probability density in space - time fractional diffusion.
% \vskip .25truecm \noindent 
%{AMS Subject Classification:}  26A33. 33E12, 45K05, 60G18,  60G50, 60G52, 60K05, 76R50.
% \vskip .25truecm \noindent
%{Keywords:} Continuous Time Random Walk, Fractional Diffusion, 
% Mittag-Leffler Function, Power Laws,   L\'evy Stable Distributions, Subordination
%%%%%%%%%%%
\end{abstract}
%%%%%%%%%%
\vsp
%% 26A33,  %%%%  (main);    Fractional derivatives and integrals
%% 33E12, %% Mittag-Leffler type functions
%% 33C60,  %% hypergeometric integrals and functions defined by them
%% 44A10,  %% Laplace Transforms
%% 45K05,  %% integro-partial differential equations
\vsp

\section*{1. Introduction}
In recent decades the notion of continuous time random walk has become popular as a  
conceptual and computational model for diffusion processes. 
Its relation to continuous diffusion processes is an asymptotic one, 
either by considering spatial and temporal refinement or by looking for large-time wide-space 
behaviour (both views being equivalent to change the units of measurement in time and space). 
Often such walks are considered for power law behaviour of waiting times and jumps, 
and from the behaviour in the Laplace-Fourier domain near zero the relevant limiting diffusion equations 
are deduced heuristically.
There are several methods of making precise this transition to the limit, 
the  fundamental role being played by the domains of attraction of stable probability laws (in time and space). 
It is desirable to have presentations of this transition without much theory of abstract probability and stochastic processes, and this can be done in a transparent manner for 
the derivation of the limiting fractional diffusion equations. 
The essential idea of our approach presented here is to use the integral equation of continuous time 
random walk separately for three random walks obtained by splitting the original walk,
 and then synthesizing the results. We so obtain our method of "parametric subordination"
as well as what usually is called "subordination" in fractional diffusion.
%%%%%%%%%%%%%%%%%
\section*{2. The Space-Time Fractional Diffusion}
 We begin by considering the Cauchy problem for the
(spatially one-dimensional) {\it space-time fractional diffusion} (STFD) equation.
$$  {\, _t}D_{*}^{\, \beta }\, u(x,t)
 \, = \,
 {\, _x}D_{ \theta}^{\,\alpha} \,u(x,t)\,,
\quad  u(x,0) = \delta (x)\,, \quad x \in \RR,\quad t \ge 0\,, \eqno(2.1) $$
where 
%% $ -\infty<x<+\infty\,,$ $\,t\ge 0\,, $
   \{$\alpha \,,\,\theta\,,\, \beta $\} are real parameters
 restricted to the ranges
$$ 0<\alpha\le 2\,,\quad  |\theta| \le \min \{\alpha, 2-\alpha\}\,,
  \quad 0<\beta\leq 1\,.\eqno(2.2) $$
Here
${ \,_t}D_*^{\,\beta}  $
denotes  the
{\it  Caputo fractional derivative}
of order $\beta $, acting on the time variable $t$,
and  $ {\,_x}D_{\,\theta}^{\,\alpha}$
denotes
the {\it  Riesz-Feller fractional derivative}
of order $\alpha $ and
skewness $\theta$,   acting on the space variable $x$.
%%%%
Let us note that the solution $u(x,t)$  of the
 Cauchy problem (2.1), known as the {\it Green function} or fundamental solution
of the space-time fractional diffusion equation,
is a
probability density in the spatial variable $x$, evolving in time $t$.
In the case $\alpha =2$ and $\beta =1$  we recover
the standard diffusion equation for which the fundamental solution
is the Gaussian density with variance $\sigma ^2 =2t$.
\vsp
Writing, with $\hbox{Re} [s] > \sigma _0$, $\kappa \in \RR$, the transforms of Laplace and Fourier as
$$     {\L} \lt\{ f(t);s\rt\}=  \widetilde f(s)
 := \int_0^{\infty} \e^{\ds \, -st}\, f(t)\, dt\,,$$
 $${\F} \lt\{g(x);\kappa\rt\}=  \widehat g(\kappa)
  := \int_{-\infty}^{+\infty} \!\! \e^{\,\ds i\kappa x}\,g(x)\, dx\,,
$$
we have the corresponding transforms
of ${ \,_t}D_*^{\,\beta} f(t)  $
and       $ {\,_x}D_{\,\theta}^{\,\alpha} g(x)$ as
$$   {\L} \lt\{  {\,_t}D_*^{\,\beta}\, f(t)\rt\} =
    s^{\, \ds \beta}\,  \widetilde{f} (s)- s^{\, \ds \beta -1}\, f(0)\,,
\eqno(2.3)$$
$$   {\F} \lt\{ {\,_x}D_{\,\theta}^{\,\alpha}\, g(x )\rt\} =
  - |\kappa |^{\,\ds \alpha}  \, i ^{\,\ds \theta \,\sgn \kappa} \,
   \widehat {g}(\kappa )\,. \eqno(2.4)$$
Notice that
$i ^{\,\ds \theta \,\sgn \kappa}= \exp [ i \,(\sgn \kappa)\,\theta\,\pi/2 ]$.
%% and $\sgn \kappa =1$ for $\kappa >0$, $=-1$ for $\kappa <0$
%% and $0$ for $\kappa =0$.
For  the mathematical details  the interested
reader is referred to
Gorenflo and Mainardi (1997), Podlubny (1999), Butzer and Westphal (2000), 
Kilbas, Srivastava and Trujillo (2006), Hilfer (2008), Tomovski, Hilfer and Srivastava (2010)
 on fractional derivatives,
and   to Samko, Kilbas and Marichev (1993), Rubin (1996), on the Feller potentials.
\vsp
For our purposes let us here confine ourselves to recall the
representation in the Laplace-Fourier domain of the (fundamental) solution
of (2.1) as it results from the application of the transforms
of Laplace and Fourier. Using $\widehat \delta (\kappa ) \equiv 1$
we have from (2.1)
 $$  s^{\, \ds \beta}\,\widehat{\widetilde{u}}(\kappa ,s) - s^{\, \ds \beta -1}
    = -|\kappa|^{\, \ds \alpha} \,  i ^{\,\ds \theta \,\sgn \kappa}
   \, \widehat{\widetilde{u}}(\kappa ,s) \,,
$$
hence
$$  \widehat{\widetilde{u}}(\kappa ,s)
    =  \frac{ s^{\, \ds \beta -1}}
{s^{\,\ds \beta} + |\kappa |^{\, \ds \alpha} \,
i^{\,\ds \theta \,\sgn \kappa} }\,.
   \eqno(2.5)$$
For explicit expressions and plots of  the fundamental solution of (2.1)
in the space-time domain
we refer the reader to Mainardi, Luchko and Pagnini (2001).
There, starting from the fact  that the Fourier transform
$\widehat{u}(\kappa ,t)$ can be written as a Mittag-Leffler function
with complex argument, the authors
have  derived a Mellin-Barnes integral representation
of $u(x,t)$  with which they have proved the non-negativity
of the solution for values of the parameters
$\{\alpha,\, \theta, \,  \beta \}$ in the range (2.2)
and analyzed the evolution in time of its  moments.
In particular for $\{0<\alpha <2, \, \beta=1\}$ we obtain
the stable densities of order $\alpha$ and skewness $\theta$.
The representation of $u(x,t)$ in terms of Fox $H$-functions
can be found  in  Mainardi, Pagnini and Saxena (2005), see also
Chapter 6 in the book by Mathai, Saxena and Haubold (2010).
%%%%%%%%%%%%%%
We note, however,  that the solution of the STFD Equation (2.1) and
its variants has been investigated by several authors;
let us only mention some of them,
 Barkai (2001),
 Meerschaert et al. (2002 and 2004), 
%%  PRE02sol,%%%
Metzler and Klafter  (2004), where the connection with the $CTRW$
was also pointed out.
%%%%%%%%CTRW
% In the present paper it is our first intention to show
% how fractional diffusion (2.1) can be obtained from continuous time
% random walks by well-scaled transition to the diffusion limit.
%\vskip -0.5truecm
\section*{3. The  Continuous-Time  Random Walk}
%%%%%%%%%%%%%%%%
Starting in the Sixties and Seventies of the past century the concept of continuous 
	  time random walk, $CTRW$, became popular in physics as a rather general (microscopic) 
	  model for diffusion processes, see the monograph by Weiss (1994). 
	  Mathematically, a $CTRW$ is 
	  a {\it renewal process with rewards} or a 
	  {\it random-walk subordinated to a renewal process}, 
	  and has been treated as such by Cox (1967).  
	  It is generated by a sequence of independent identically distributed ($iid$)
	   positive waiting times $T_1,T_2,T_3, \dots$, 
	  each having the same probability distribution $\Phi(t)$, $t\ge0$, 
	  of waiting times $T$, and a sequence of $iid$ 
	  random jumps $X_1, X_2, X_3, \dots$, 
	  each having the same probability distribution function $W(x)$, $x \in  \RR$,
	  of jumps $X$. 
	  These distributions are assumed to be independent of each other. 
	  Allowing generalized functions (interpretable as measures), 
	  %% in the sense of Gel`fand and Shilov (1964) 
	  we have corresponding probability densities
	  $\phi(t)=\Phi'(t)$  and $w(x)= W'(x)$ that we will use for ease of notation. 
	  Setting $t_0=0$, $t_n=T_1+T_2+\dots + T_n$  for $n\in \NN$, and
	  $x_0=0$, $x_n= X_1+X_2+\dots + X_n$, $x(t)=x_n$  
	  for $t_n\le t<t_{n+1}$,  we get a (microscopic) model of a diffusion process. 
	  A wandering particle starts in  $x_0=0$ and makes a jump $X_n$ {\it at each instant} $t_n$.
	  %%%%%%%%%%%%%%% 
	  Natural probabilistic reasoning then leads us to the {\it integral equation 
	  of continuous time random walk}  for the probability density $p(x,t)$  of 
	  the particle being in position $x$  at instant $t\ge 0$:
$$
 p(x,t) =  \delta (x)\, \Psi(t)\, +
  \int_0^t  \!\!  \phi(t-t') \, \left[
 \int_{-\infty}^{+\infty}\!\!  w(x-x')\, p(x',t')\, dx'\right]\,dt'\,,
\eqno(3.1)$$                
Here the {\it survival probability} 
$$\Psi(t) = \int_t^\infty \phi(t') \, dt'
\eqno(3.2)$$                                          
denotes the probability that at instant $t$ the particle still is sitting 
in its initial position $x_0=0$. 
Eq. (3.1) is also known as (semi --) Markov renewal equation.
Germano et al. (2009), considering the more general situation of coupled CTRW, 
have given a detailed stochastic proof of an integral equation which by modification and 
specialization reduces to (3.1).
\vsp
By using the Fourier and Laplace transforms
we arrive, via   $\widehat \delta(\kappa) \equiv 1$
and the convolution  theorems, in the transform domain at the equation
$$
\widehat{\widetilde p}(\kappa ,s)
 =   \widetilde\Psi(s) +
   \widehat w(\kappa )\,\widetilde \phi(s) \,
   \widehat{\widetilde p}(\kappa ,s) \,, 
\eqno(3.3)$$        
which, by $\widetilde \Psi(s)= (1-\widetilde \phi(s))/s$
    implies the Montroll-Weiss equation, see e.g  Weiss (1994),
$$\widehat{\widetilde p}(\kappa ,s) = 
\frac{1-\widetilde \phi(s))}{s} \, \frac{1}{1- \widehat w(\kappa )\,\widetilde \phi(s)}\,. 
\eqno(3.4)$$                       
Because of $|\widehat w(\kappa)| <1$,
$|\widetilde \phi(s)|<1$   for  $\kappa \ne 0$, $s \ne 0$,  we can expand 
%% into a geometric series
$$
\widehat{\widetilde p}(\kappa ,s) = \widetilde\Psi(s)\, \sum_{n=0}^{\infty}
   \,\left[\widetilde \phi(s) \,\widehat w(\kappa )\right]^n\,,
\eqno(3.5)$$                 
and promptly obtain, %% the series representation of the $CTRW$, 
see e.g.  Cox (1967) and  Weiss (1994),
$$
p(x,t)= \sum_{n=0}^{\infty} v_n(t)\, w_n(x)\, \q -\infty<x<+\infty, t\ge 0\,.
\eqno(3.6)$$                     
Here the functions  $v_n = \Psi\,* \, \phi^{*n}$ and $w_n= w^{*n}$ 
are obtained by iterated convolutions in time $t$ and in space $x$, respectively; 
in particular we have: 
$$ v_0(t) =(\Psi \,*\, \delta)(t) = \Psi(t), \;
v_1(t) =(\Psi \,*\, \phi)(t) , \; w_0(x)=\delta(x),\; w_1(x)= w(x).\eqno(3.7)$$
The representation (3.6) can be found without the detour over (3.4) by direct probabilistic treatment. 
It exhibits the $CTRW$ as a subordination of a random walk to a renewal process.    
  \vsp
   Note that in the special case $\phi(t)= m \,\exp(-mt)$, $m>0$,
    the equation (3.1) describes the compound Poisson process. 
	It reduces after some manipulations %% (best carried out in the transform domain) to
	to the {\it Kolmogorov-Feller equation}
     $$ \frac{\d}{\dt} p(x,t)= -m p(x,t) + m \int_{-\infty}^{+\infty} w(x-x')p(x',t)\, dx'\,. \eqno(3.8)$$                         
From (3.6)  we then obtain the well-known series representation
      $$  p(x,t)=\e^{-mt}\, \sum_{n=0}^{\infty}\frac{(mt)^n}{n!} \, w_n(x)\,.\eqno(3.9)$$  

\section*{4. Power Law Assumptions and Well-Scaled Passage to the Diffusion Limit}
We will pass to the diffusion limit of $CTRW$ under power law assumptions
 for waiting times and jumps given by asymptotic relations %% in the Laplace - Fourier domains
$$ 1 - \widetilde {\phi} (s) \sim \lambda s^\beta ,\quad
  \lambda >0,  \quad  s \to 0^+\,, \eqno(4.1)$$
$$ 1-\widehat{w}(\kappa ) \sim \mu
|\kappa| ^\alpha \, i^{\ds \, \theta \,\sgn \kappa }\,,
 \quad  \mu >0\,, \quad  \kappa \to 0\,,
 \eqno(4.2) $$
where $\beta$, $\alpha$ and $\theta$ are restricted as in (2.2).
See, e.g. Gorenflo and Mainardi (2003),  Scalas, Gorenflo and Mainardi (2004), Gorenflo (2010) for their physical
  meaning, 
obtainable via Tauber-Karamata theorems on aymptotics.
%% `With these assumptions 
We rescale space and time by multiplying the jumps  $X_k$ and waiting times $T_k$
by positive factors $h$ and $\tau$ (intended to be sent to zero), so obtaining a rescaled walk
$$ x_n (h)= \sum_{k=1}^n h \, X_k\,, \q t_n(\tau)=\sum_{k=1}^n \tau\, T_k\,.$$
Physically, this amounts to changing the units of measurement from $1$ to $1/h$ in space $x$,
from $1$ to $1/\tau$ in time $t$.
\vsp
For the reduced jumps and waiting times we get
$$  w_h(x)= w(x/h)/h\,, \q \phi_\tau (t)= \phi(t/\tau)/\tau\,, \eqno(4.3)$$
$$ \widehat w_h(\kappa)= \widehat w(h\kappa),, \q \widetilde \phi_\tau (s) = \widetilde \phi(\tau s)\,, \eqno(4.4)$$
and the Montroll-Weiss equation (3.3) transforms into
$$
\widehat{\widetilde p}_{h,\tau}(\kappa ,s) \!=\! 
\frac{1-\widetilde \phi_\tau (s)}{s} \, \frac{1}{1- \widehat w_h(\kappa )\,\widetilde  \phi _\tau(s)}
\!=\! \frac{1-\widetilde \phi (\tau s)}{s} \, \frac{1}{1- \widehat w(h\kappa )\,\widetilde \phi(\tau s)}\,.
\eqno(4.5) $$
Using our power law assumptions (4.1), (4.2) we find, for $h\to 0$, $\tau \to 0$, the asymptotics
$$
\widehat{\widetilde p}_{h,\tau}(\kappa ,s) \sim 
{\ds \frac{\lambda \tau^\beta s^{\beta-1}}{\mu (h |\kappa|)^\alpha\, 
i^{\ds \theta\, \sgn \kappa} + \lambda (\tau s)^\beta}}
= {\ds \frac{s^{\beta-1}}{\rho(h,\tau)\, |\kappa|^\alpha\, i^{\ds \theta\, \sgn \kappa} + s^\beta}}
\eqno(4.6) $$  
with
$$ \rho(h,\tau)  = \frac{\mu h^\alpha}{\lambda \tau^\beta}\,.$$
Adopting the scaling relation
$$ \rho(h,\tau)  = \frac{\mu h^\alpha}{\lambda \tau^\beta} = 1\,,\eqno(4.7)$$
we arrive for $h \to 0$ at the limit
 $$ \widehat{\widetilde p}_{0,0}(\kappa ,s) = 
	 {\ds \frac{s^{\beta-1}}{|\kappa|^\alpha \, i^{\ds \theta\, \sgn \kappa}+ s^\beta}} =:
	\widehat{\widetilde u}(\kappa,s)   $$
	which we identify as the Fourier-Laplace solution to the Cauchy problem (2.1) for the STFD
	 equation.
	By inverting the Laplace transform and then for the Fourier transform invoking the continuity theorem of probability theory
	we get weak convergence to the fundamental solution of (2.1). 
%%%%%%%%%%\newpage
%%%%%%%%
\vsp 
{\bf Remark}.
 Because Mittag-Leffler waiting times distributed according to %%  (see (3.2))
$\Psi(t) = E_\beta (-t^\beta)$ play a distinguished role among power law waiting times and 
lead to time--fractional CTRW, in fact to the time--fractional Kolmogorov--Feller equation, 
see e.g. Hilfer (1995), Hilfer and Anton (1995), Gorenflo and Mainardi (2008) and Gorenflo (2010),
 it is natural to use them (properly scaled) for approximately simulating particle paths 
 of space-time fractional diffusion processes. 
 Fulger, Scalas and Germano  (2008) have done this,
  thereby using exact methods for generating the required temporal Mittag-Leffler deviates and spatial 
  stable deviates, giving references for these methods.        
  %%%%%%%%%%
\section*{5. Splitting into Three Random Walks}
We analyze the anatomy of the $CTRW$ of Section 3 by splitting it into three random walks 
($rw_1$), ($rw_2$), ($rw_3$)
 that we will treat in analogy to our method of well-scaled transition to the diffusion limit.
 We will work in three coordinates: the semiaxis $t\ge 0$ of natural time,
 the semi-axis $t_*\ge 0$ of operational time,
 the axis $-\infty<x<+\infty$ of space.
 In addition to the Laplace transform $\L\{f(t);s\} = \widetilde f(s)$
 with respect to $t$ we will need the Laplace transform $\L\{f(t_*);s_*\} = \widetilde f(s_*)$
 with respect to $t_*$.
 \vsp
 For convenience let us denote by ($rw$) $\,\equiv\,$ ($CTRW$) the random walk of Section 3,
 $$ t_0=0\,,\;  t_n = T_1 + T_2 +  \dots + T_n\,; 
 \q  x_0 =0\,, \; x_n = X_1 + X_2 + \dots  + X_n \,.\eqno(5.1)$$
   This is a random walk along the space line $x\in \RR$ happening in physical time $t\ge 0$.
   Here $\phi(t)$ is the waiting time density, $w(x)$ is the jump density, 
   and $p(x,t)$ is the sojourn density in point $x$ at instant $t$. 
   
   By ($rw_1$) we denote the basic renewal process
   $$ t_{*,n}=n\, (n= 0, 1, 2, \dots)\,; 
   \q      t_0 =0\,,\; t_n = T_1 + T_2 +  \dots + T_n\,. \eqno(5.2)$$
   We consider it as a random walk along natural time $t\ge 0$, happening in operational 
   time $t_*\ge 0$. Here $t_*$ plays the role of a ``pseudo time", $t$   
 plays the role of ``pseudo space". 
 The waiting time (in $t_*$) is \underline{constant} (fixed) =1, so we have the
``waiting time" density $v(t_*)= \delta(t_*-1)$. 
 The jump instants are the integers $t_{*n}= n>0$,
and we set $t_{*0}=0$.
%%%%%%%%%% 
 The ``jump density" (in $t$) is $\phi(t)$. 
 We denote the ``sojourn density"  in ``point" $t$ at ``instant" $t_*$
 by  $r(t, t_*)$.
 Note that in this  sojourn density $r$ the first argument is the pseudo space $t$, 
 whereas the second argument is pseudo time $t_*$.
 \vsp
 By ($rw_2$) we denote the random walk  
 $$
 t_{*,n}=n\, (n= 0, 1, 2, ...)\,; 
 \q x_0=0\,,\; x_n = X_1 + X_2 + \dots + X_n \,,\eqno(5.3)$$
 which we consider as a random walk along natural space 
 $x\in \RR$, happening in operational time $t_*\ge 0$ which (again) plays the role of "pseudo time". 
 But $x$ here is the usual spatial variable.
 The waiting time density  in $t_*$ is (again) $v(t_*)= \delta(t_*-1)$,
  the jump density is $w(x)$,    and we denote the sojourn density in point $x$ at "instant" $t_*$ by
 $\xi(x,t_*)$.
 \vsp
 By ($rw_3$) we denote the random walk
 $$ t_0=0\,,\;  t_n = T_1 + T_2 +  \dots + T_n\,;
  \q t_{*n} =n \, (n= 0, 1, 2, ...)\,,\eqno(5.4)$$
 along the semi-axis $t_*\ge 0$ happening  in physical time $t\ge 0$.
 The walk is just the inversion of ($rw_1$). 
 The variable $t$ now plays the role of time, the waiting time density is $\phi(t)$.
 But $t_*$ plays the role of "pseudo space", the jump density being $v(t_*)= \delta(t_*-1)$.
 The sojourn  density in "point" $t_*$ at instant $t$ is denoted 
 by $q(t_*,t)$ where the pseudo space $t_*$ is the first argument and time $t$ is the second argument. 
 
\section*{6. Asymptotic Analysis of the Three Walks}
In ($rw_1$) and ($rw_3$) we have walks in the positive pseudo-space direction.
 For such walks it is convenient to modify the Montroll-Weiss formula by using 
 {\it the Laplace transform in place of the Fourier transform}.
 If in the general ($CTRW$) the jump density $w(x) \ge 0$ for all $x\in \RR$ then with
 $$ \widetilde g(\kappa):= \int_0^\infty \!\! g(x)\, \e^{-\kappa x}\, dx\,, \q
 \kappa \ge \kappa_0\,,\eqno(6.1)$$
 we can replace the original Montroll-Weiss equation (3.4) by
$$\widetilde{\widetilde p}(\kappa ,s) = 
\frac{1-\widetilde \phi(s)}{s} \, \frac{1}{1- \widetilde w(\kappa )\,\widetilde \phi(s)}\,. 
\eqno(6.2)$$        
Here for ($rw_1$) and ($rw_3$) the specifically relevant replacements
of $x$ and $t$, $\kappa$ and $s$ must be made, analogously for ($rw_2$)  in (3.3) and (3.4).
Analogous to Section 4 we will carry out well-scaled transitions to the diffusion limit, observing
(4.1)  and (4.2) and an additional asymptotic relation for the Laplace transform 
of $v(t_*)= \delta(t_* -1)$.
We have
$$ \widetilde v(s_*) = \e^{-s_*} = 1 - s_* + o(s_*)\,, \q \hbox{for}\q s_*\to 0\,.$$
%%%%%%%%%%%%%%%%%%%%%%%%%%%%%%%%
\subsection*{6.1 Analysis of ($rw_1$).}
%%%%%%%%%%%%%%%%%%%%%%%%%%%%%%%%%%
In place of (6.2) we have
$$ \widetilde{\widetilde r}(s ,s_*) = 
\frac{1-\widetilde v(s_*)}{s_*} \, \frac{1}{1- \widetilde \phi(s)\,\widetilde v(s_*)}\,. 
\eqno(6.3)$$
We rescale in $t_*$ direction by a factor $\sigma>0$, so replacing the constant waiting time 1
by $\sigma$. We replace the ``jumps" $T$ by $\tau T$ with $\tau>0$.    
So we obtain, for $\tau \to 0$, $\sigma \to 0$:
$$  
\begin{array}{ll}
   {\ds \widetilde{\widetilde r}_{\tau,\sigma}(s ,s_*)} & = 
{\ds \frac{1-\widetilde v(\sigma s_*)}{s_*} \, \frac{1}{1- \widetilde \phi(\tau s)\,\widetilde v(\sigma s_*)}}\\
& \sim  {\ds \frac{\sigma}{\lambda \tau^\beta s^\beta + \sigma s_*}}= 
{\ds \frac{1}{\frac{\lambda \tau^\beta}{\sigma}\, s^\beta + s_*}}\,.
\end{array}
$$
By introducing the scaling relation
$$ \frac{\lambda \tau^\beta}{\sigma} = 1 \,,\eqno(6.4)$$
we arrive for $\tau \to 0$ at the limit
$$  \widetilde{\widetilde r}_{0,0}(s ,s_*) = \frac{1}{s^\beta + s_*}\,, $$
whence
$$\widetilde r_{0,0}(s ,t_*) = \e^{-t_*s^\beta}\,, $$
corresponding to 
$$ r_{0,0}(t ,t_*) = t_*^{-1/\beta}\, L_\beta^{-\beta}\left(t_*^{-1/\beta}\, t\right)\,,
\eqno(6.5)$$
which is a probability density in $t$ 
   evolving with $t_*$.
Here $L_\beta^{-\beta}(t)$ denotes the right-extremal stable density whose Laplace transform is $\exp(-s^\beta)$.
%%%%%%%%%%%%%%%%%%%%%%%%%%%%%%%%
\subsection*{6.2 Analysis of ($rw_2$).}
%%%%%%%%%%%%%%%%%%%%%%%%%%%%%%%%%%
In analogy to (3.4) we find
$$
\widehat{\widetilde \xi}(\kappa ,s_*) \!=\! 
\frac{1-\widetilde v (s_*)}{s_*} \, \frac{1}{1- \widehat w(\kappa )\,\widetilde  v(s_*)}\,,
\eqno(6.6)
 $$
and by rescaling (replacing ``waiting times" 1 by $\sigma>0$, jumps $X$ by $hX$)
$$
\begin{array}{ll}
{\ds \widehat{\widetilde \xi}_{h,\sigma}(\kappa ,s_*)} &=
{\ds \frac{1-\widetilde v (\sigma s_*)}{s_*} \, \frac{1}{1- \widehat w(h\kappa )\,\widetilde  v(\sigma s_*)}}\\
&{\ds \sim \frac{\sigma}{\mu h^\alpha |\kappa|^\alpha i^{\ds \theta \, \sgn \kappa} +\sigma s_*}}
 = {\ds \frac{1}{\frac{\mu h^\alpha}{\sigma}\, |\kappa|^\alpha i^{\ds \theta \, \sgn \kappa} +\sigma s_*}}\,.
 \end{array}$$
 The scaling relation needed here is
 $$ \frac{\mu h^\alpha}{\sigma}=1\,,    \eqno(6.7)$$
 and with it we get the limit
$$ \widehat{\widetilde \xi}_{0,0}(\kappa ,s_*) =
 \frac{1}{s_* +|\kappa|^\alpha i^{\ds \theta \, \sgn \kappa}}\,,$$
 $$\widehat \xi_{0,0}(\kappa ,t_*) =  \exp \left( -t_* |\kappa|^\alpha \,i^{\ds \theta \, \sgn \kappa}\right)\,,$$
 hence 
 $$ \xi_{0,0}(x,t_*) = t_*^{-1/\alpha}\, L_\alpha^{\theta}\left(t_*^{-1/\alpha}\, x\right)
 = f_{\alpha,\theta}(x, t_*)
\,, \eqno(6.8) $$
which is a probability density in $x$ evolving with $t_*$.
Here $L_\alpha^{\theta}(x)$ is the stable density with Fourier transform
$\exp \left( -|\kappa|^\alpha \,i^{\ds \theta \, \sgn \kappa}\right)$
and  $f_{\alpha,\theta}(x,t)$ is the  solution $u(x,t)$ of the Cauchy problem
       (2.1) for $\beta = 1$.
%%%%%%%%%%%%%%%%%%%%%%%%%%%%%%%%

%%%%%%%%%%%%%%%%%%%%%%%%%%%%%%%%
\subsection*{6.3 Analysis of ($rw_3$).}
%%%%%%%%%%%%%%%%%%%%%%%%%%%%%%%%%%
In place of (6.2) we get
$$\widetilde{\widetilde q}(s_* ,s) = 
\frac{1-\widetilde \phi(s)}{s} \, \frac{1}{1- \widetilde v(s_*)\,\widetilde \phi(s)}\,. 
\eqno(6.9)$$
 Replacement of waiting times $T$ by $\tau T$, "jumps" 1 by $\sigma$ yields
 $$
\begin{array}{ll}
{\ds \widetilde{\widetilde q}_{\sigma,\tau}(s_*,s)} &=
{\ds \frac{1-\widetilde \phi (\tau s)}{s} \, \frac{1}{1- \widetilde v(\sigma s_*)\,\widetilde \phi(\tau s)}}\\
&{\ds \sim \frac{\lambda \tau^\beta s^{\beta-1}}{\sigma s_* +\lambda \tau^\beta s^\beta}}
 = {\ds \frac{s^{\beta-1}}{\frac{\sigma}{\lambda \tau^\beta} s_* + s^\beta}}\,.
 \end{array}
 \eqno(6.10)$$
 Using the same scaling relation as for ($rw_1$), namely $\lambda \tau^\beta/\sigma =1$, we get
 the limit
 $$ \widetilde{\widetilde q}_{0,0}(s_*,s) = \frac{s^{\beta-1}}{s_* + s^\beta}\,, \eqno(6.11)$$
 implying 
 $$ \widetilde q_{0,0}(t_*,s) = s^{\beta-1}\, \e^{-t_* s^\beta}\,, \eqno(6.12)$$
 so that
 $$  q_{0,0}(t_*,t) = \, _tJ^{1-\beta}\, r_{0,0}(t, t_*)\,,\eqno(6.13)$$
 with $r_{0,0}$ as in (6.5) and the operator $\,_tJ^{\gamma}$, $\gamma>0$
 of Riemann-Liouville fractional integration from $0$ to $t$ whose Laplace symbol is $s^{-\gamma}$.
 Note that  $q_{0,0}(t_*,t)$ is a probability
   density in $t_*$ evolving with $t$.
\vsn   
 {\bf Remark}. We also recognize that   $q_{0,0}(t_*,t)$ can be expressed in terms of 
 the so-called $M$-Wright function, see e.g. the recent book by Mainardi (2010), Appendix F, 
 Eqs (F.13), (F.24) and (F.51)--(F.53).
  In fact, inverting (6.12), we get, in view of (F.51) and (F.52),
 $$  q_{0,0}(t_*,t)= t^{-\beta}\, M_\beta(t_*/t^\beta)\,.\eqno (6.14)$$
 To the same result we arrive inverting (6.11) from $s$ to $t$, getting
 $$ \widetilde q_{0,0}(s^*,t) = E_\beta(-s^*t^\beta)\,, \eqno(6.15)$$
 where $E_\beta$ is the Mittag-Leffler function. Then, inverting this Laplace transform from $s^*$ to $t^*$
 we get, in view of (F.53), the result (6.14).

 %%%%%%%%%%%%%%%%%
\section*{7. Synthesis for Subordination}
%%%%%%%%%%%%%%%%%%%%%%%%
The scaling relations $\sigma =\lambda \tau^\beta$ for ($rw_1$) and ($rw_3$) and $\sigma= \mu h^\alpha$
 for (rw2) imply $ \mu h^\alpha = \lambda \tau^\beta$,
 the scaling relation (4.7) for well-scaled transition to the diffusion limit of ($CTRW$).  
%% \vsp
Synthetizing the limiting processes $y =y(t_*)$ of ($rw_2$) and $t_*=t_*(t)$ of ($rw_3$) we get the limiting
process $x=x(t)$ of ($CTRW$) as $x=y(t_*(t))$, which implies the integral formula of
 subordination for the solution of the 
Cauchy problem (2.1), namely
$$ u(x,t) = \int_{t_*=0}^\infty \!\!f_{\alpha, \theta}(x, t_*)\, q(t_*,t)\, dt_*\,, \eqno(7.1)$$
with
$$q(t_*,t)= q_{0,0}(t_*,t) = t_*^{-1/\beta}\,_tJ^{1- \beta}\, L_\beta^{-\beta}\left(t_*^{-1/\beta}\, t\right)
= t^{-\beta}\, M_\beta(t_*/t^\beta) \,,
\eqno (7.2)$$
which is a probability density in $ t_ *$ evolving with $t$, see (6.13)-(6.14).
%%%%%%%%%%%%%%%%%%%%
%%  \vsp %%%%%%%%%
%% the following sentence has been inserte by RG on September 18.
Note that we denote the parent process by $y(t_*)$ in distinction to the limit $x(t)$ of ($CTRW$). 
Compare with Meerschaert et al. (2002) and with Gorenflo, Mainardi and Vivoli  (2007).
\vsp
The limit of ($rw_1$) and ($rw_2$) yields our {\it method of parametric subordination}. 
They describe the processes $t=t(t_*)$ and $y=y(t_*)$, hence, by identifying $y$ with the spatial variable 
$x$ and eliminating $t_*$, our method of parametric subordination 
$x= y(t_*)$, $t = t(t_*)$ for constructing trajectories of particles.  
In this method we call $y=y(t_*)$ the ``{\it parent process}", 
$t=t(t_*)$ the ``{\it leading process}".
Leaning on Feller  (1971)  we call $t_*=t_*(t)$ the ``{\it directing process}".
\vsp
Let us remark that the methods described by Fogedby (1994), Kleinhans and Friedrich (2007),
Zhang, Meerschaert and Baeumer (2008)  working with Langevin equations 
%% and circumventing the CTRW concept 
are, in a certain sense,  equivalent to parametric subordination.
In the two coupled stochastic differential equations (6.3) and (6.4) 
 of Gorenflo, Mainardi and Vivoli (2007) we have hinted to this 
 Langevin approach without calling it so.
 Fogedby (1994) describes in detail this approach for the spatially 
 multi-dimensional situation and discusses its relevance for treatment 
 of problems from applications. Kleinhans and Friedrich (2007) produce and  
 exhibit several particle trajectories and discuss information that can 
 be drawn from these. Zhang, Meerschaert and Baeumer (2008) do likewise and 
 consider furthermore situations of diffusion type and with 
 variable coefficients in the spatial operator. As a concrete 
 application they outline what is happening in Bromide plumes occurring in 
 a site of North America, actual irregular processes that have inspired 
 hydrological research. They compare results of measurements with those 
 of their simulations and so demonstrate convincingly the power of 
 their theory.
\vsp
  It is our intention to compare our approach to simulation of particle 
 trajectories  with  other approaches  and generalizations in a 
 forthcoming paper.  The essence of our approach is that it is based on  
 a systematic and consequent application of the CTRW integral equation 
 to the various processes involved.

\vskip -0.5truecm 

%%%%%%%%%%%%%%%%%%%%%%%%%
%\paragraph{Notes by Gorenflo}

% Fig. 0 shows in form of a diagram the connections between the four random walks (rw1),
% (rw2),( rw3) and (CTRW)=(rw). We can consider Fig. 1 as a representation of (rw2),
% Fig. 2 as one of (rw1) or by interchange of axes as one of (rw3), and finally Fig. 3
% as a representation of (CTRW).

% Unfortunately, the enumeration of the random walks is different from that of the
% figures. But we should leave it as such. There is not enough time to change it, and
% such change is dangerous, it can introduce new errors.

% In our paper with Vivoli $t=t(t_*)$ is called the "leading process" in contrast to
% the "directing process" $t_* = t_*(t)$. But, differing from our paper with Vivoli and
% also from our Section 7 on Synthesis, in the present paper in Fig. 2 the process
% $t=t(t_*)$ is now called the "directing" process. For the sake of consistency I would
% be happy if this could be changed. The inverse process $t_*=t_*(t)$ is the directing
% process, also in Feller's terminology (difficult to see and easy to confuse due to
% pure mathematicians like Feller and Samko and many others using the same letter for
% different things). If the text to Fig. 2 cannot be changed in short time, then we
% should just leave it as it is. 
% Probably nobody considers the matters so carefully. 
% We can then be extremely careful in our paper for Mathai.

\section*{8. Numerical Results}
%%%%%%%%%%%%%%%%\vfill\eject%%%%%%%%%%%%
%% The method of producing the stable
%% random deviates should be described or at least a source should be
%% given from where it is taken.
%% Vivoli has adopted
%% the method by Janicki, see
%% \cite{Janicki_LN96,Janicki-Weron_94}
%%%%%%%%%%%
In this  Section,
after describing the numerical schemes adopted,
we shall show the sample paths for two case studies
of symmetric ($\theta=0$) fractional diffusion processes:
%% $\{\alpha =2,  \, \beta =0.90\}$,
$\{\alpha =2,  \, \beta =0.80\}$,
$\{\alpha =1.5, \,  \beta =0.90\}$.
%% $\{\alpha =1.5, \,  \beta =0.80\}$.
As explained in the previous sections, for each case
we need to construct the sample paths for three  distinct processes,
the leading process $t = t(t_*)$ ,
the parent process $x= y(t_*)$ (both in the operational time)
and, finally, the subordinated process $x =x(t)$,
corresponding to the required fractional diffusion process.
We shall depict the above sample paths in Figs. 1,  2,  3,
respectively, devoting the left  and the right plates to
the different case studies.
%% In order to produce fractional diffusion sample paths
For this purpose, following Gorenflo, Mainardi and Vivoli (2007), we  proceed as follows.
\vsp
First, let the operational time $t_*$ assume only integer
values, say  $t_{*,n}=n$ with $n=0,1,\dots,10000$.
Then, produce $10000$ independent
identically distributed random deviates,
say $X_1,X_2,\dots,X_{10000}$, having a symmetric
stable probability  distribution of order $\alpha$.
%% and skewness $\theta,$
% see the book by Janicki \cite{Janicki_LN96}
% for a useful and efficient method to do that.
Now, with the  points
 $$x_0=0\,,\;  x_n = \sum \limits _{k=1}^n X_k,\quad n\ge 1\,,\eqno(8.1)$$
the couples $(t_{*,n},x_n)$, plotted
in the  $(t_*,x)$ plane  (operational time, physical space)
can be considered as points of
a sample path
$\{y(t_*) : 0\le t_*\le 10000\}$
of a symmetric L\'evy motion with order  $\alpha$
corresponding
to the integer values of operational time $t_*=t_{*,n}$.
In this identification of $t_*$ with $n$ we use the fact that
our stable laws for waiting times and jumps imply
$\lambda = \mu = 1$  in the asymptotics (4.1) and (4.2) and
$\tau = h = 1$ as initial scaling factors in (4.3).
\vsp
In order to complete the sample path we agree to connect every
two successive  points
$(t_{*,n},x_n)$  and  $(t_{*,n+1},x_{n+1})$
by a horizontal line from $(t_{*,n},x_n)$ to
$(t_{*,n+1},x_n),$ and a vertical line from
$(t_{*,n+1},x_n)$ to $(t_{*,n+1},x_{n+1}).$
%%%%%%%%%%%%
Obviously, that is not the ``true" L\'evy motion
from point  $(t_{*,n},x_n)$ to point
$(t_{*,n+1},x_{n+1}),$
but from the theory of $CTRW$ we know this kind of sample
path to converge to the corresponding  L\'evy motion paths in the diffusion limit.
%%%%%%%%% INSERTION by GORENFLO, 16 May 05
However, as the successive values of $t_{*,n}$ and $x_n$ are
generated by successively adding the relevant standardized stable
random deviates,  the obtained sets of points in the three coordinate
planes: $(t_*, t)$, $(t_*,x)$, $(t,x)$ can, in view of  infinite
divisibility and self-similarity of the stable probability distributions,
be considered as
snapshots of the corresponding true random processes occurring in
continuous operational time $t_*$ and physical time $t$,
correspondingly.
Clearly, fine details between successive points are missing.
\vsp
The well-scaled passage to
the diffusion limit here consists simply in regularly subdividing the
$\{t_*\}$ intervals  of length 1   into smaller and smaller
 subintervals (all of equal length) and adjusting the random
increments of $t$ and $x$ according to the requirement of  self-similarity.
%%%%%%%%%%%%
Furthermore if we watch
a sample path in a large interval of operational time $t_*,$
the points $(t_{*,n},x_n)$
and  $(t_{*,n+1},x_{n+1})$ will in the graphs appear very near
to each other  in operational time $t_*$
and aside
from missing mutually cancelling jumps up and down (extremely near to
each other) we have a good picture of the true processes.
%%%%%%%%%%%%%%%%%%%%%%
%%\vsp
\begin{figure}
\includegraphics[width=.48\textwidth]{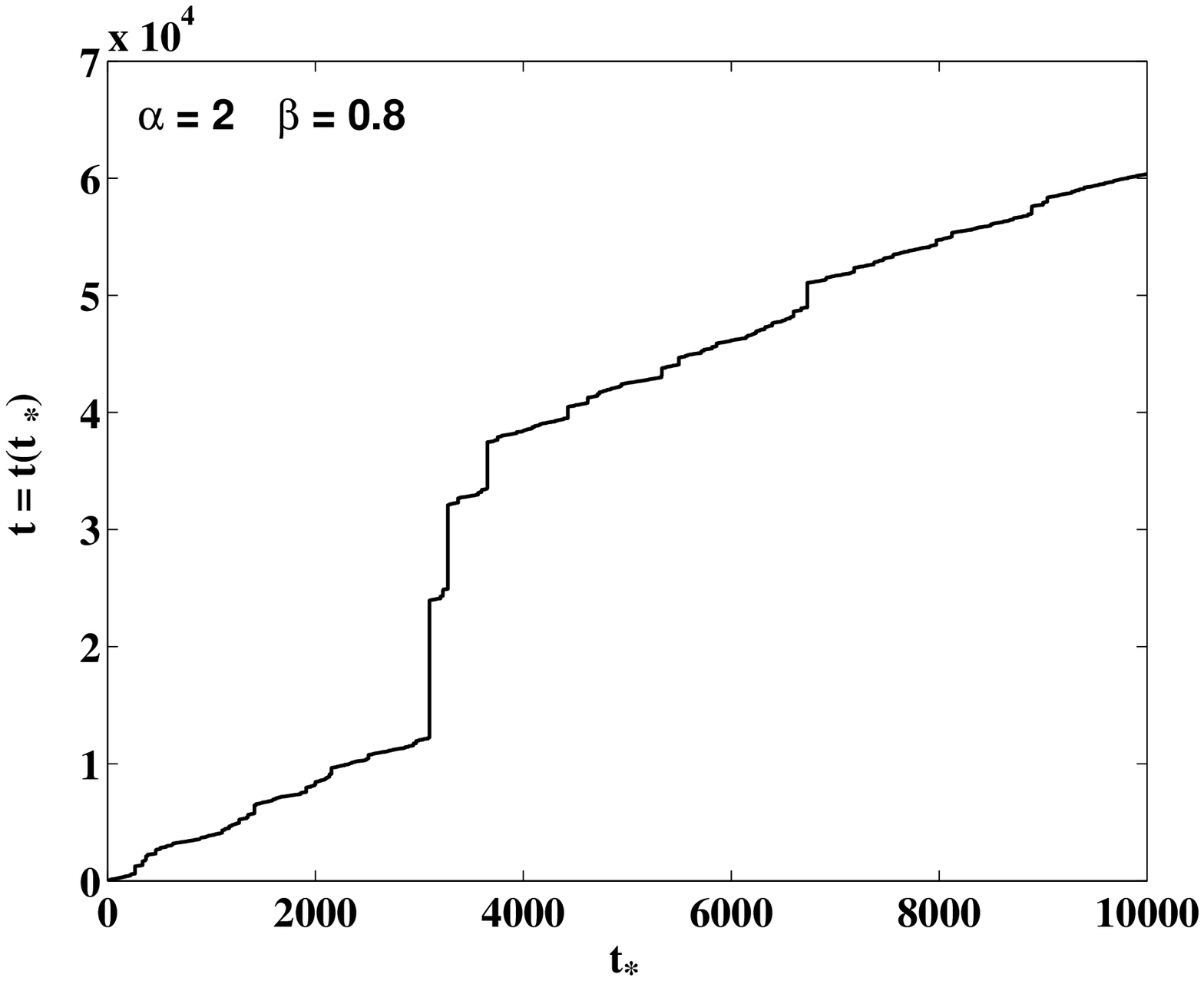}
 \includegraphics[width=.48\textwidth]{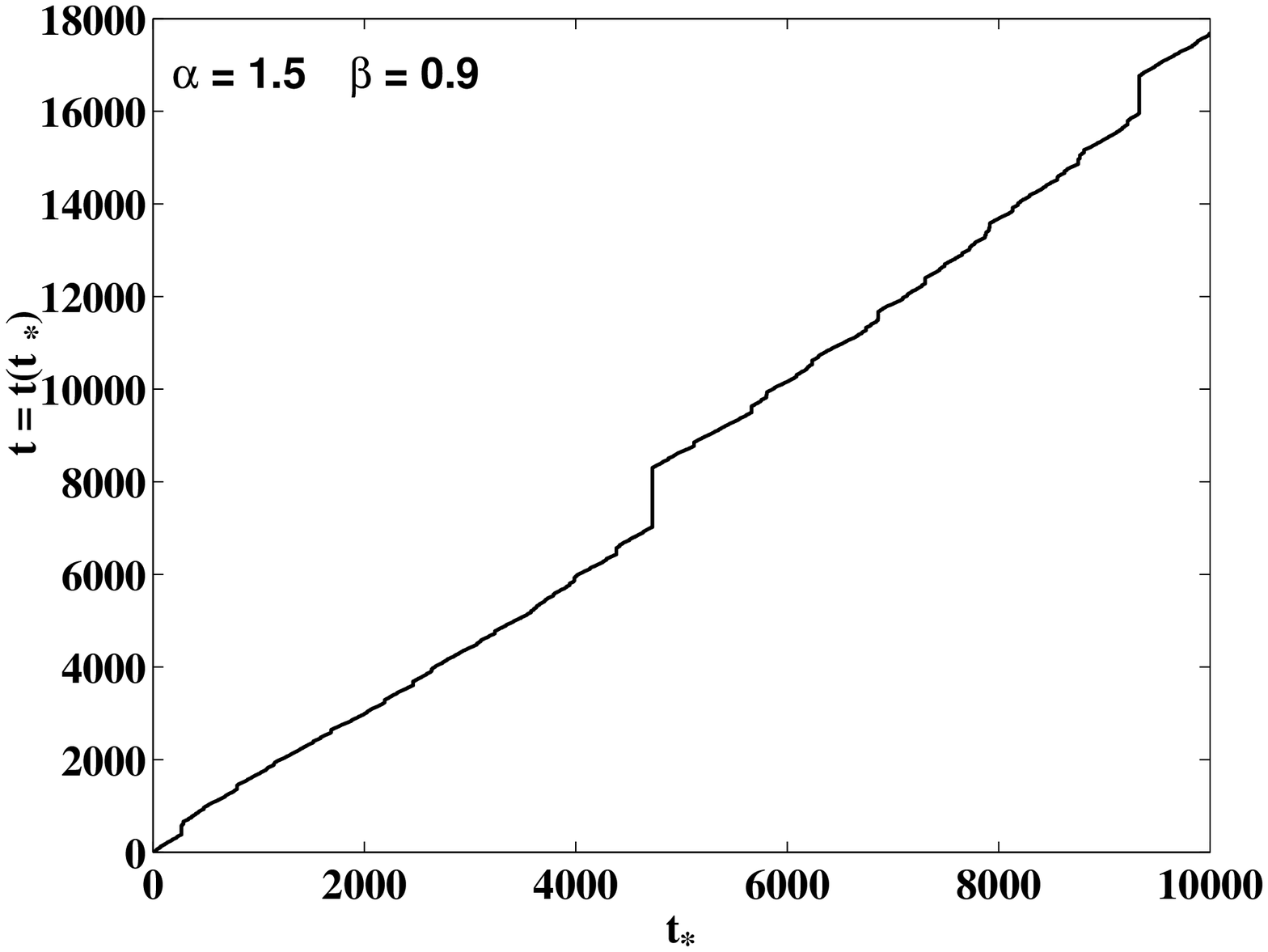}
 \caption{A sample path for ($rw_1$), the leading process $t=t(t_*)$.}
 \centerline{LEFT: $\{\alpha =2\,,\; \beta =0.80 \}$,
      RIGHT: $\{\alpha =1.5\,,\; \beta =0.90 \}$.}
%%%\end{figure}
\vskip 0.75truecm
 \includegraphics[width=.48\textwidth]{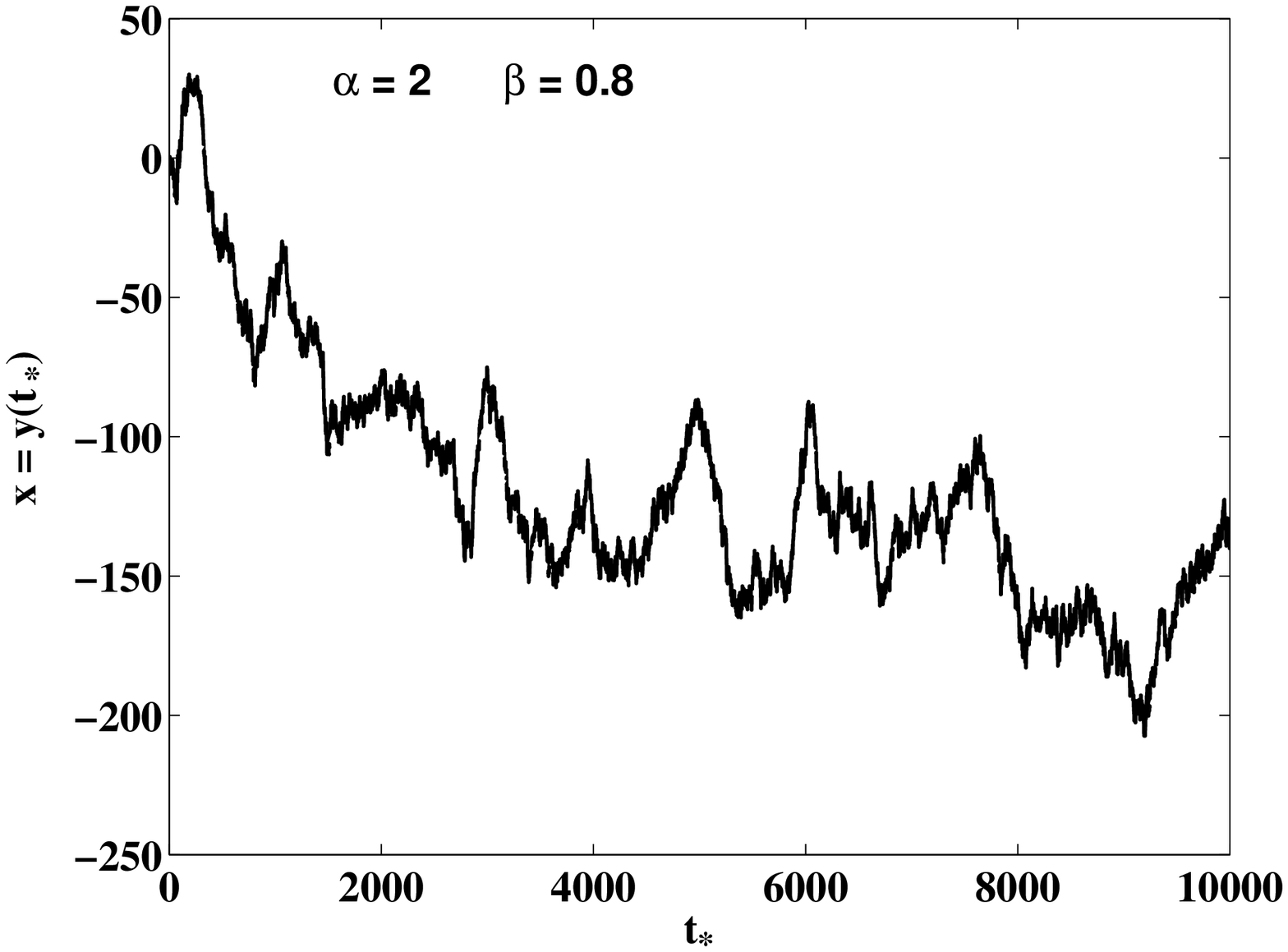}
 \includegraphics[width=.48\textwidth]{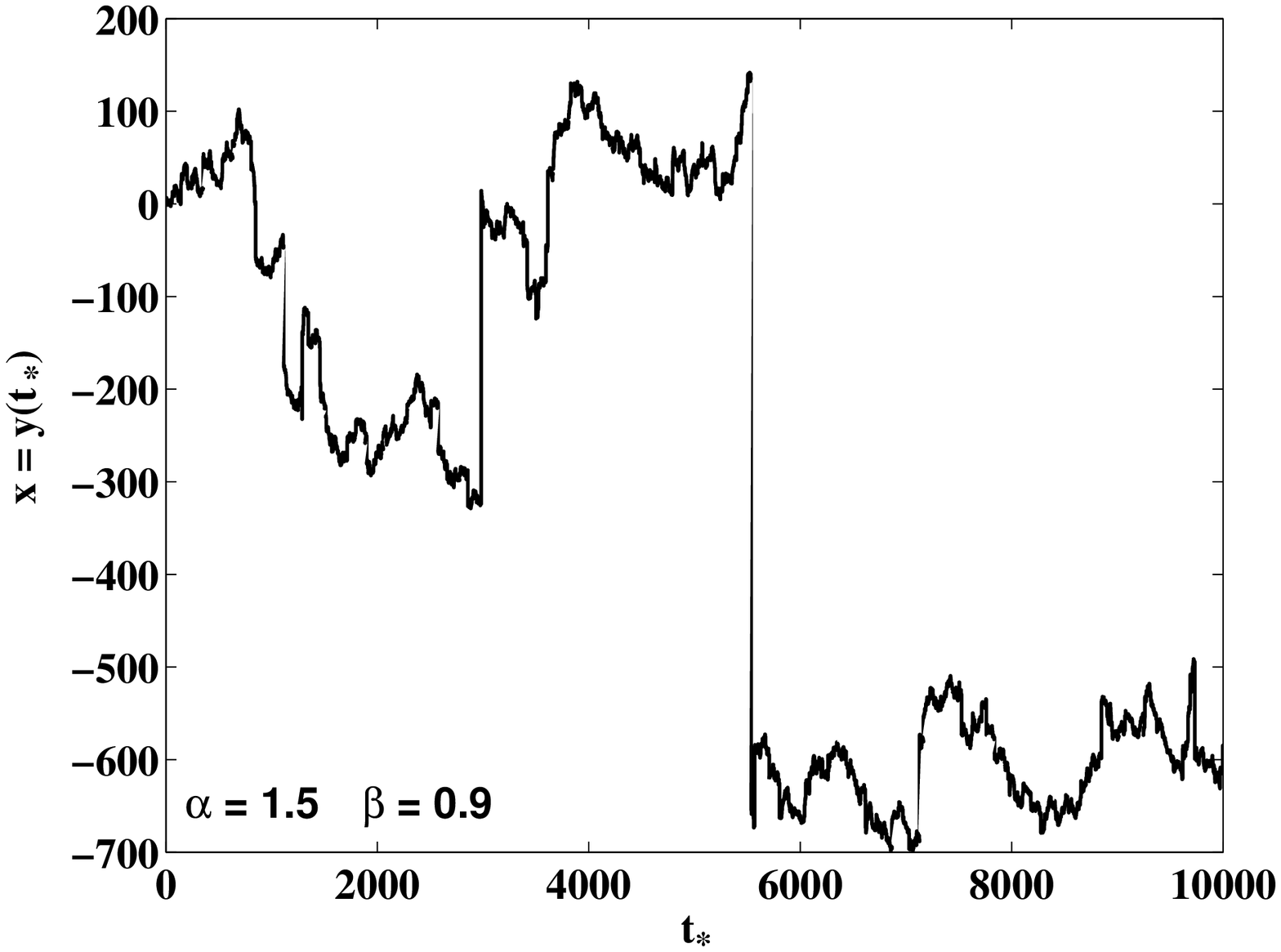}
 \caption{A sample path for ($rw_2$), the parent process $x=y(t_*)$.}
 \centerline{LEFT: $\{\alpha =2\,,\; \beta =0.80 \}$,
      RIGHT: $\{\alpha =1.5\,,\; \beta =0.90 \}$.}
%%%\end{figure}
%%%\begin{figure}
\vskip 0.75truecm
\includegraphics[width=.48\textwidth]{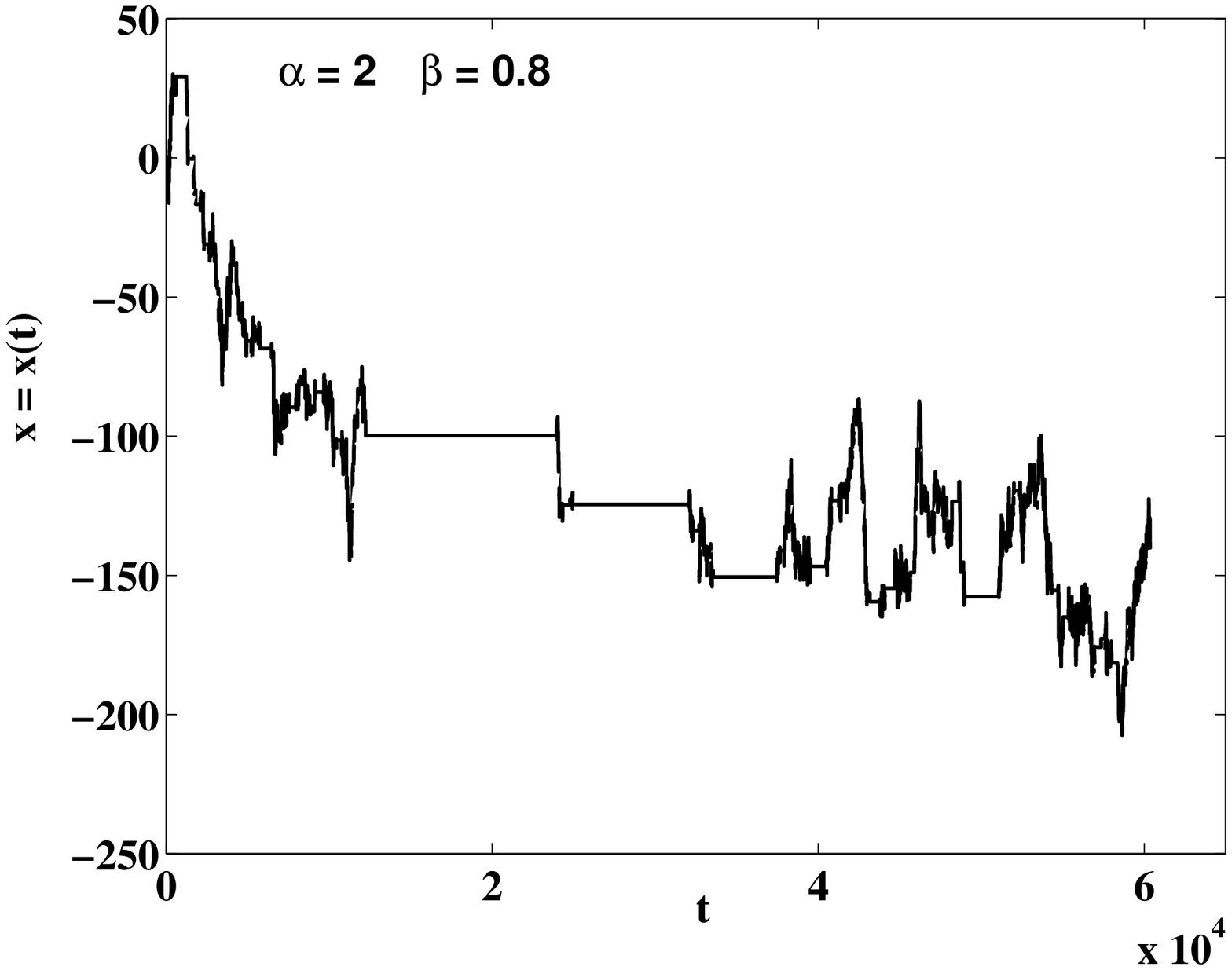}
 \includegraphics[width=.48\textwidth]{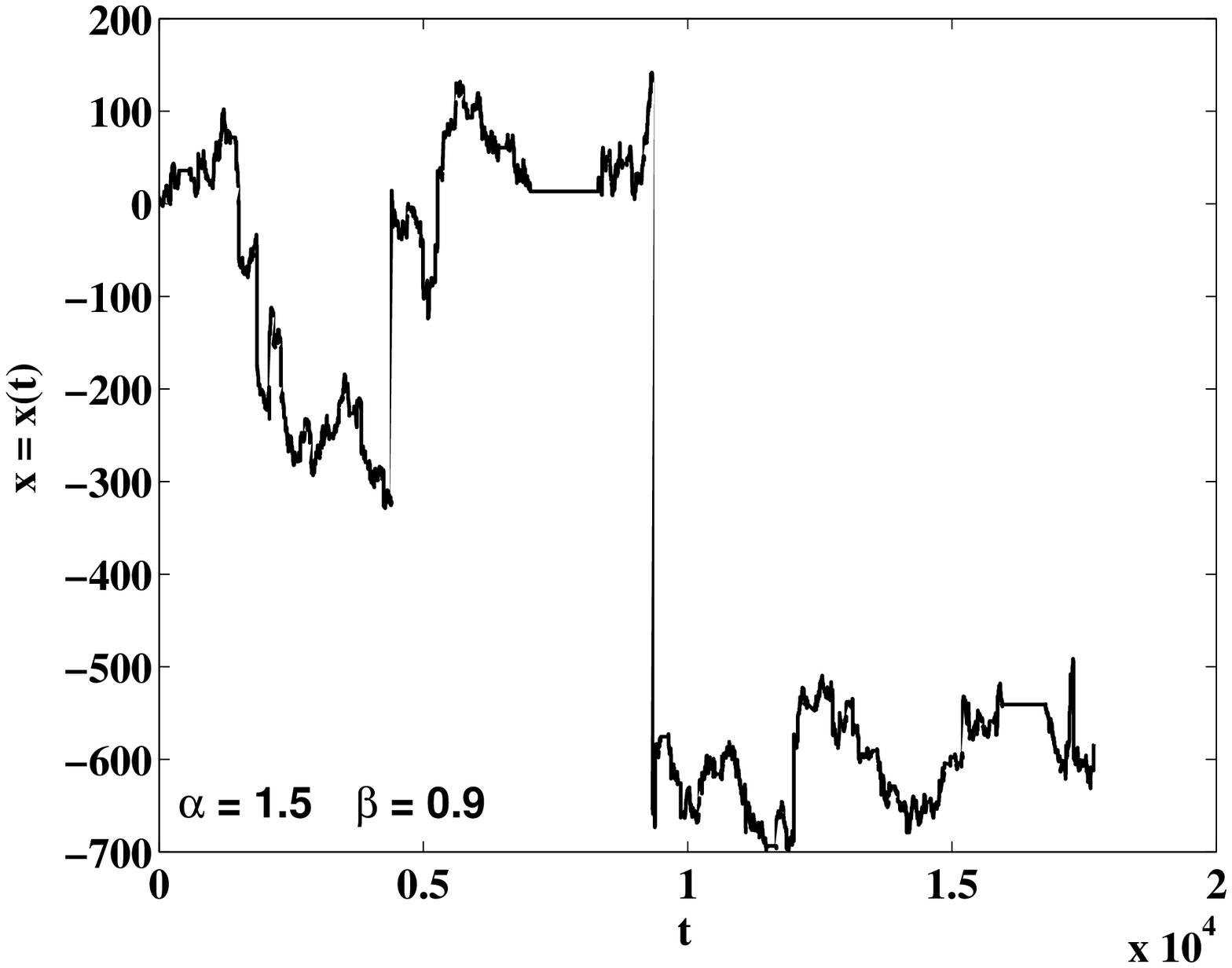}
 \caption{A sample path for ($rw$), the subordinated process $x=x(t)$.}
 \centerline{LEFT: $\{\alpha =2\,,\; \beta =0.80 \}$,
      RIGHT: $\{\alpha =1.5\,,\; \beta =0.90 \}$.}
\end{figure}
%%%%%%%%%%%%%
\noindent
%%%%%%%%%%
Plots in Fig. 1 (devoted to the leading process, the limit of ($rw_1$))
thus represent sample paths in the $(t_*,t)$  plane
of unilateral L\'evy motions  of order $\beta$. 
By interchanging the coordinate axes we can consider Fig 1 as representing sample paths of 
the directing process, the limit of ($rw_3$). 
\vsp
Plots in Fig. 2  (devoted to the parent process, the limit of ($rw_2$))
represent   sample
paths in the $(t_*,x)$ plane,
produced in the way explained above, for  L\'evy motions
of  order  $\alpha$ and skewness  $\theta=0$
(symmetric stable distributions).
%% keeping fixed parameter $\theta=0.$
% \vsp
As indicated above, it is possible to produce $10000$ independent
identically distributed random deviates,
say $T_1,T_2,\dots,T_{10000}$ having a stable
probability  distribution with order $\beta$
and skewness $ -\beta$ (extremal stable distributions).
Then, consider the points
 $$t_0=0,\quad  t_n = \sum\limits_{k=1}^n T_k,\quad n\ge 1\,,\eqno(8.2)
$$
and plot
the couples $(t_{*,n},t_n)$
in the  $(t_*,t)$ (operational time, physical time) plane.
By  connecting  points with
horizontal and vertical lines we get   sample paths
$\{t(t_*) : 0\le t_*\le 10000\}$
%% with order $\beta$ and skewness $-\beta,$
describing the evolution
of the physical time $t$ with the increasing of the operational time $t_*.$
%%%%%%%%%
%\vsp
%%%%%%%%%%%%%%%%%%%%%
Now, plotting  points $(t(t_{*,n}),y(t_{*,n}))$ in the  $(t,x)$ plane,
namely the physical time-space plane, and connecting them as before,
one gets a good
approximation of
the sample  paths of the subordinated fractional diffusion
process of parameters
$\alpha$, $\beta$ and $\theta=0$.
%%%%%%%%%%%%%%%%%
\vsp
In Fig. 3 (devoted to the subordinated process, the limit of ($CTRW$))
we show paths obtained in this way, from
the points calculated in the previous paths.
In Fig.3 Left we plotted the points $(t(t_{*,n}),y(t_{*,n}))$
obtained in Fig. 1 Left and Fig. 2 Left,
while Fig.3 Right shows the points of Fig. 1 Right and Fig. 2 Right.
\vsp
By observing the figures the reader will note that
horizontal segments (waiting times) in the $(t, x)$ plane (Fig. 3)
correspond to vertical segments (jumps) in the $(t, t_*)$ plane (Fig. 1).
%%%%%%%   Comments added by Gorenflo 30-31 October 2005
Actually, the graphs in the $(t, x)$-plane depict continuous time random
walks with waiting times $T_k$ (shown as horizontal segments) and jumps
$X_k$ (shown as vertical segments). The left endpoints of the horizontal
segments can be considered as snapshots of the true particle path (the true
random process to be simulated), the segments being segments of our ignorance.
In the interval $t_n < t \le t_{n+1}$ the true process
(namely the spatial variable $x=X(t)$)
may jump up and down (infinitely) often, the sum (or integral) of all these
ups (counted positive) and downs (counted negative) amounting to the
vertical jump $X_{n+1}$.
%%%%%%%%%%
%% \vsp
Finer details will become visible by choosing in the operational time
$t_*$ the step length $\tau  << 1$ (instead of length 1 as we have done) and
 correspondingly the  waiting times and spatial jumps as $\tau^{1/\beta}$
multiplied by a standard extreme $\beta$-stable deviate, $\tau^{1/\alpha}$
multiplied by a standard (in our special case: symmetric) $\alpha$-stable
deviate, respectively, as required by the self-similarity properties of
the stable probability distributions.
%%%%%%%%%%\vsp
In a forthcoming paper we will describe in detail what happens for finer
and finer discretization of the operational time $t_*$.

\vskip -0.5truecm

\section*{9. Conclusions}
%%%%%%%%%%%%%%%
\begin{figure}[h!]
\vspace{-0.5truecm}
\begin{center}
\includegraphics[width=.60\textwidth]{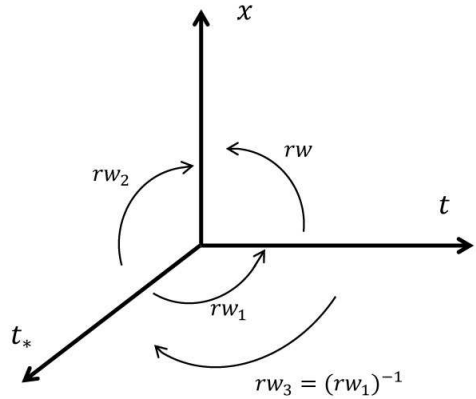}
\end{center}
\vspace{-1truecm}
 \caption{Diagram for the connections between the four random walks ($rw_1$),
($rw_2$),($rw_3$) and ($CTRW$)=($rw$). }
\end{figure}
\noindent
After sketching our method of well-scaled transition to the diffusion limit under power law regime,
we have slpitted the $CTRW$ into three random walks($rw_1$), ($rw_2$), 
($rw_3$), each of them involving two of the three directions:
natural time, operational time and space. The random walk ($rw_3$) is the inverse
of ($rw_1$).
To better visualize the matter,   here  Fig. 4 shows in form of a diagram the connections 
between  the three random walks ($rw_1$), ($rw_2$), ($rw_3$) and ($rw$) $\,\equiv\,$ ($CTRW$).
%%  where we have denoted by ($rw_3$) the inverse of $rw_1$.
\newpage
\vsp
In section 8 the figures 1, 2, 3 have shown some numerical realizations of these random walks.  
We can consider Fig. 1 as a representation of ($rw_1$),
or by interchange of axes as one of ($rw_3$),
Fig. 2 as one of ($rw_2$), and finally Fig. 3
as a representation of ($CTRW$).
\vsp
We have carried out the analysis of transition by aid of the transforms of Fourier and Laplace,
thereby consistently using corresponding variants of the Montroll-Weiss equation 
(not only in space and natural time, but also not conventionally in the above-mentioned combinations).
In this way we have exhibited the true power of this equation in asymptotic analysis of random walks.
In a quite natural way we have interpreted the diffusion limits as leading to the common subordination formula 
as well as the method of parametric subordination of constructing trajectories of particles.
%%%%%%%%%%%%%
\vskip -0.5truecm
\subsection*{Acknowledgments} 
We gratefully appreciate helpful comments of our colleague Enrico Scalas to a draft of this paper.

%%%%%%%%%%%%%%%%
\vskip -0.5truecm

%%%%%%%%%%%%%%%%%%%%%%
% \begin{thebibliography}{99}
\section*{References}

\vsn
 1. E. Barkai,
 {Fractional Fokker-Planck equation, solution, and application},
 {\it Phys. Rev. E} {\bf 63}, 046118-1/18 (2001).

\vsn   %% \bibitem[Butzer and Westphal (2000)]{Butzer-Westphal 00}
2. P. Butzer, U. Westphal,
Introduction to fractional calculus,
 in: H. Hilfer (Editor),
{\it Fractional Calculus, Applications in Physics}
 (World Scientific, Singapore, 2000), pp. 1--85.

\vsn
3. W. Feller, %%  (1971).
{\it An Introduction to Probability Theory and its Applications}, Vol II
  (Wiley, New York, 1971).
  %% Second Edition. [First edition (1966)]

 \vsn
4.  H.C. Fogedby, 
 Langevin equations for continuous time L\'evy flights,
 {\it Phys. Rev. E} {\bf 50}, 1657--1660 (1994).
 
 \vsn
5. D. Fulger, E. Scalas, G. Germano,  
 Monte Carlo simulation of uncoupled continuous-time random walks yielding a 
 stochastic solution of the space-time fractional diffusion equation. 
{\it Phys. Rev. E} {\bf 77}, 021122/1--7 (2008).

\vsn
6. G. Germano, M. Politi, E. Scalas, R.L. Schilling, 
 Stochastic calculus for uncoupled continuous-time random walks. 
 {\it Phys. Rev. E} {\bf 79}, 066102/1--12 (2009).

\vsn
7.  R. Gorenflo, Mittag-Leffler waiting time, power laws, rarefaction,
continuous time random walk, diffusion limit, in     
 S.S. Pai, N. Sebastian, S.S. Nair, D.P. Joseph and D. Kumar (Editors),
%% CMS Proceedings Pala July 2010, you have been sent this volume.
Proceedings of the National Workshop on  Fractional Calculus and Statistical Distributions,
%% Edited by Shanoja S. Pai, Nicy Sebastian, Seema S. Nair, Dhannya P. Joseph, Dilip Kumar
%% http://www.cmsintl.org/general_announcements/WorkshopOnFractionalCalculus_Proceedings.pdf
%% 25-27 November 2009,
(CMS Pala Campus, India, 2010), pp. 1--22.
[E-print: {\tt http://arxiv.org/abs/1004.4413}] %% , 28 pages.
%  \bibitem{GorMai_CISM97}  %% [3]
 \vsn
8.  R. Gorenflo, F. Mainardi,  
  Fractional calculus: integral and differential equations of fractional order,
  in: A. Carpinteri and F. Mai\-nardi (Editors),
  {\em Fractals and Fractional Calculus in Continuum Mechanics\/}
  (Springer Verlag, Wien  1997),  pp. 223--276. %% [CISM Lecture Notes Vol. 378]
  [E-print: {\tt http://arxiv.org/abs/0805.3823}]
%%  [Reprinted in NEWS 010101, see {\tt http://www.fracalmo.org}]
%%   

% \bibitem{GorMai_INDIA03}    %% [4]
 \vsn
9. R. Gorenflo,  F. Mainardi,
 Fractional diffusion processes: probability distributions and continuous time random walk, in:
    G. Rangarajan and M. Ding (Editors),
    {\it Processes with Long Range Correlations}
    (Springer-Verlag, Berlin, 2003), pp. 148--166. %% [Lecture Notes in Physics, No. 621]
   [E-print: {\tt http://arxiv.org/abs/0709.3990}]
\vsn
10. R. Gorenflo, F. Mainardi, %%  (2008): 
Continuous time random walk, Mittag-Leffler waiting  time and fractional diffusion: 
mathematical aspects,   in: 
R. Klages, G.  Radons,  and I.M. Sokolov,  (Editors),  
{\it Anomalous  Transport, Foundations and Applications} 
 (Wiley-VCH Verlag, Weinheim, Germany, 2008),  pp. 93--127. %%  ISBN 978-3-527-40722-4.
 [E-print: {\tt arXiv:cond-mat/07050797}]
 %% Paper presented at the WE-Heraeus-Seminar on Anomalous Transport: 
%% Experimental  Results and Theoretical Challenges, Physikzentrum Bad-Honnef (Germany),  
%% 12-16 July 2006.

%% \bibitem{GorMaiViv_CSF07} %% [5]
\vsn
11. R. Gorenflo, F. Mainardi, A. Vivoli,
Continuous time random walk and parametric subordination in fractional diffusion,
   {\it Chaos, Solitons and Fractals} {\bf 34},  87--103 (2007).
   [E-print http://arxiv.org/abs/cond-mat/0701126]

\vsn
12. R. Hilfer, Threefold introduction to fractional calculus, 
 %% Paper presented at the WE-Heraeus-Seminar on Anomalous Transport:
  %% Experimental Results and Theoretical Challenges, Physikzentrum
%% Bad-Honnef (Germany), 12-16 July 2006.  Published on pages 17-73 in:
%% Rainer Klages, Günter Radons, and Igor M. Sokolov (editors):
in R. Klages, G. Radons and I.M. Sokolov (Editors),
{\it Anomalous  Transport, Foundations and Applications}
(Wiley-VCH Verlag, Weinheim, Germany, 2008), pp. 17--73.
\vsn
13.   %%  \bibitem{Hilfer_FRACTALS95}
 R. Hilfer,
  Exact solutions for a class of fractal time random walks,
 {\it  Fractals} {\bf 3}, 211-216 (1995).
%%%%%%%%%%%%%%%%%%%%%%
 \vsn
14. %%  \bibitem{Hilfer-Anton_PRE95}
 {R. Hilfer, L. Anton},
 Fractional master equations and fractal time random walks,
 {\it Phys. Rev. E} {\bf 51}, R848--R851 (1995).
\vsn  %% \bibitem{Kilbas-Srivastava-Trujillo_BOOK06}  %% [6]
15. A.A. Kilbas, H.M. Srivastava, J.J. Trujillo, 
{\it Theory and Applications of Fractional Differential Equations}
(Elsevier, Amsterdam, 2006). 
% [North-Holland Series on Mathematics Studies No 204]
%%%%%%%%%
\vsn
16. D. Kleinhans, R. Friedrich, 
 Continuous-time random walks: Simulations of continuous trajectories, 
 {\it Phys. Rev E} {\bf 76}, 061102/1--6 (2007).

\vsn %%% \bibitem{Mainardi BOOK10}
17. F. Mainardi, 
{\it Fractional Calculus and Waves in Linear Viscoelasticity}
 (Imperial College Press, London, 2010).

% \bibitem{Mainardi_FCAA01}  %% [7]
\vsn
18. F. Mainardi, Yu. Luchko, G. Pagnini,
     The fundamental solution of the space-time fractional diffusion
     equation,
   {\it Fract. Calculus and Appl. Analysis} {\bf 4}, 153--192 (2001).
 [E-print: {\tt http://arxiv.org/abs/cond-mat/0702419}]
%% [Reprinted in NEWS 010401, see {\tt http://www.fracalmo.org}]
% [Paper dedicated to Professor Rudolf Gorenflo for his 70-th birthday]

%% \bibitem{Mainardi_JCAM05}   %% [8]
\vsn
19. F. Mainardi, G. Pagnini, R.K. Saxena,
 Fox $H$ functions in fractional diffusion,
 {\it   J. Computational and Appl. Mathematics} {\bf 178},  321--331 (2005).
%% Proceedings of the 7-th International Symposium on
%% Orthogonal Polynomials, Special Functions and Applications (OPSFA),
%% Copenhagen (DK) 18-22 August 2003.  URL: www.math.ku.dk/conf/OPSFA2003
%%%%%%%%%%%%%%%%
\vsn %% \bibitem{Mathai-Saxena-Haubold BOOK-H-2010}
20. A.M Mathai, R.K. Saxena, H.J Haubold,
{\it The H-function, Theory and Applications}
(Springer Verlag, New York, 2010).
%%%%%%%%

% \bibitem{M3_PRE02sub}  %% [9]
\vsn
21. M.M. Meerschaert, D.A. Benson, H.P  Scheffler, B. Baeumer,
Stochastic solutions of space-fractional diffusion equation,
{\it Phys. Rev. E} {\bf  65},  041103/1--4 (2002).
%%A PRESENTATION, GOOD AND SHORT,FOR SUBORDINATION
%%%%%%%%%%%%%%%%%%
% \bibitem{M3_PRE02sol}  %% [10]
\vsn
 %% M.M. Meerschaert, D.A. Benson, H.P  Scheffler and P. Becker-Kern,
 %% Governing equations and solutions of anomalous random walk limits,
 %% {\it Phys. Rev. E} {\bf  66},  060102-1/4 (2002).
%%%%%%%%%%%%%%%%%%%%%%
22. M.M. Meerschaert, H.P  Scheffler,
Limit theorems for  continuous-time random walks with infinite mean waiting times, 
{\it J. Appl. Prob.} {\bf  41}, 623--638 (2004).
%%%%%%%%%%%%%%%%%%%%%
\vsn %%  \bibitem{Metzler-Klafter_JPhysics04}  %% [11]
23. R. Metzler, J. Klafter,
The restaurant at the end of the random walk: Recent developments
 in the description of anomalous transport by fractional dynamics,
 {\it J. Phys. A. Math. Gen.}  {\bf 37},  R161--R208 (2004).

\vsn %% \bibitem{Podlubny_99}   %% [12]
24.  I. Podlubny,
  {\it Fractional Differential Equations}
  (Academic Press, San Diego, 1999).
%%%%%%%%%%%%%
% \bibitem{Saichev_PhysA05}
% A. Piryatinska, A.I. Saichev and W.A. Woyczynski,
% Models of anomalous diffusion: the subdiffusive case,
%  {\it Physica A}  {\bf 349},  375-420 (2005).
%%%%%%%%%%%%%%%%%%%

\vsn %% \bibitem[Rubin (1996)]{Rubin BOOK96}
25.  B. Rubin, 
 {\it Fractional Integrals and Potentials}
 (Addison-Wesley \& Longman, Harlow, 1996).
%% [Pitman Monographs and  Surveys in Pure and Applied Mathematics No 82]

\vsn %% \bibitem{SKM_93}  %5 [13]
26. S.G.  Samko, A.A. Kilbas, O.I. Marichev,
{\it Fractional Integrals and Derivatives: Theory  and  Applications}
(Gordon and Breach, New York, 1993).
% Translation from the Russian edition.
%% Nauka i Tekhnika, Minsk (1987).
%%%%%
%%%%%%%%%%%%%%%%%%
% \bibitem{SGM_00}
% E. Scalas, R. Gorenflo and F. Mainardi,
% Fractional calculus and continuous-time finance,
% Physica A {\bf 284},  376-384 (2000).
%%%%%%%%%%%%%%%%%%%%%%%%%%%%%%%%%%%%%
\vsn %% \bibitem{Scalas_PRE04} %% [14]
27. E. Scalas, R. Gorenflo, F. Mainardi,
    Uncoupled continuous-time random walks:
Solution and limiting behavior of the master equation,
{\it Phys. Rev. E} {\bf 69},  011107/1--8 (2004).
%% (30 January 2004)
%%%%%%

\vsn %% \bibitem{Tomovski-Hilfer-Srivastava ITSF10}
28. Z. Tomovski, R. Hilfer, H.M. Srivastava,
Fractional and operational calculus with generalized fractional
derivative operators and Mittag-Leffler type functions,
{\it Integral Transforms Spec. Funct.} {\bf 21}, 797-- 814 (2010).
 %% No 11

\vsn %% \bibitem{Weiss_BOOK94} %% [15]
29. G.H. Weiss,
 {\it Aspects and Applications of Random Walks}
 (North-Holland, Amsterdam, 1994).
%%%%%%%%
% \bibitem{West_BOOK03}
% B.J. West, M. Bologna and P. Grigolini.
% {\it Physics of Fractal Operators}, Springer Verlag, New York (2003).
%[Institute for Nonlinear Science]
%%%%%
\vsn
30. Y. Zhang, M.M. Meerschaert, B. Baeumer,
Particle tracking for time-fractional diffusion,
{\it Phys. Rev. E.} {\bf 78}, 036705/1--7 (2008).
%%%%%%%%%%

%% \end{thebibliography}

\end{document}